\documentstyle[leqno]{article}
\newtheorem{th}{Theorem}[section]

\newcounter{defin}[section]
\renewcommand{\thedefin}{\thesection.\arabic{defin}}
\newcounter{ex}[section]
\renewcommand{\theex}{\thesection.\arabic{ex}}
\newcounter{rem}[section]
\renewcommand{\therem}{\thesection.\arabic{rem}}
\newcommand{\R}{I\!\!R}
\newcommand{\E}{I\!\!E}
\newcommand{\T}{I\!\!\!T\!\!\!\!I}
\title{CANONICAL NONLINEAR CONNECTIONS\\
ON JET BUNDLES OF FIRST ORDER}
\author{Mircea Neagu}
\date{}
\begin{document}
\maketitle
\begin{abstract}
The aim of this paper is to open the problem of finding of a nonlinear
connection $\Gamma=(M^{(i)}_{(\alpha)\beta}\;,N^{(i)}_{(\alpha)j})$
on 1-jet bundle $J^1(T,M)$, which to be canonically produced
from a given Kronecker $h$-regular fundamental vertical metrical d-tensor
$
G^{(\alpha)(\beta)}_{(i)(j)}=h^{\alpha\beta}g_{ij}\;,
$
possibly provided by multi-time dependent quadratic Lagrangians coming
from various branches of theoretical physics: bosonic strings theory
\cite{Got-Mars}, magneto-hydrodynamics \cite{Mars-Pek}, electrodynamics
\cite{Mir-An} or elasticity \cite{Olv2}.
From geometrical point of view, the importance of this problem comes
from contravariant Riemann-Lagrange geometry of 1-jet spaces \cite{Ne8},
and consists in the possibility of construction of
distinguished 1-forms
$
\delta x^i_\alpha=dx^i_\alpha+M^{(i)}_{(\alpha)\beta}dt^\alpha+N^{(i)}_
{(\alpha)j}dx^j,
$
necessary for simple local descriptions of geometrical or physical objects
studied.
\end{abstract}

\section{Geometrical and physical aspects}

\hspace{5mm} From a physical point of view, we underline that
the jet fibre bundle of order one $J^1(T,M)$ appears as a natural house
for geometrical studies of important physical domains (continuum mechanics
\cite{Mars-Pek},
quantum field  theories \cite{Olv1}, generalized multi-time relativity and
electromagnetism \cite{Ne9} or dynamical relativistic multi-time
optics \cite{Ne10}) including natural processes characterized by
dependence on position, multi-momentum or partial velocities. At the same
time, geometrical studies of first order quadratic Lagrangians from several
important branches of theoretical physics
were required a profound analysis of the differential geometry of 1-jet
spaces, in the sense of connections, torsions and curvatures.
In this direction, extending by Riemannian methodes geometrical
results from theory of Lagrange spaces \cite{Mir-An},
the paper \cite{Ne8} developed recently a contravariant geometry on
$J^1(T,M)$, in the sense of d-connections, torsion and curvature d-tensors,
which is unlike covariant Hamiltonian formalism \cite{Got-Mars},
\cite{Mars-Pek} in many ways. This geometry is naturally produced by a given
fundamental vertical metrical d-tensor
having the Kronecker $h$-regular product form

\begin{equation}\label{fdt}
G^{(\alpha)(\beta)}_{(i)(j)}(t^\gamma,x^k,x^k_\gamma)=h^{\alpha\beta}
(t^\gamma)g_{ij}(t^\gamma,x^k,x^k_\gamma),
\end{equation}
together with a nonlinear connection
$\Gamma=(M^{(i)}_{(\alpha)\beta}, N^{(i)}_{(\alpha)j})$ on jet bundle
$J^1(T,M)$.
\medskip\\
\addtocounter{rem}{1}
{\bf Remarks \therem} i) The preceding Kronecker $h$-regular vertical metrical
d-tensor may be provided by particular multi-time dependent quadratic
Lagrangian functions $L$ on 1-jet bundle $J^1(T,M)$, via the
formula
$
G^{(\alpha)(\beta)}_{(i)(j)}=(1/2)(\partial^2L/\partial x^i_\alpha\partial
x^j_\beta)
$.

ii) In our opinion, the metrical d-tensor (\ref{fdt})
may be viewed as an {\it unified} gravitational field on $J^1(T,M)$,
realized by a {\it temporal} gravitational field $h_{\alpha\beta}(t^\gamma)$
and a {\it spatial} gravitational field $g_{ij}(t^\gamma,x^k,x^k_\gamma)$
depending on multi-moment, position and partial velocities $x^i_\alpha$.
Moreover, the $x^i_\alpha$-dependence may be combined with the abstract
concept of {\it partial directions anisotropy}.\medskip

Concerning the novelty brought in
theoretical-physics, we consider that the
main and unpublished feature of Riemann-Lagrange geometry on $J^1(T,M)$
is the construction of a large geometrical background for a generalized
multi-time field theory, in the sense of generalized Maxwell and Einstein
equations, that allows the including of famous equations of
mathematical-physics (classical Maxwell or Einstein equations) as particular
cases.
We recall that the construction of a new field theory, described in
multi-time terms on $J^1(T,M)$, was required by geometrical studies of certain
famous relativistic invariant equations involving many time variables (chiral
fields, sine-Gordon), and of $KP$-hierarchy of integrable systems in
which the arbitrary variables $t^\alpha$ and $t^\beta$
are quite equal in rights and there is no reason to prefer one to another by
choosing it as time \cite{Dickey}.

From geometrical point of view, a central role in Riemann-Lagrange
differential geometry is played by nonlinear connections
$\Gamma=(M^{(i)}_{(\alpha)\beta}, N^{(i)}_{(\alpha)j})$ on $J^1(T,M)$,
that allow the construction of {\it adapted bases} of vector or covector
fields \cite{Ne8}.\medskip\\
\addtocounter{rem}{1}
{\bf Remark \therem} Concerning the physical aspects of nonlinear
connections on 1-jet spaces, we believe that these prescribe {\it possible
intrinsic interactions} between the temporal and spatial gravitational fields.
\medskip

In our opinion, the importance of adapted bases attached to nonlinear connections on
$J^1(T,M)$ may be justified in various ways.

Firstly, the simple tensorial local transformations of their elements
\cite{Ne1} imply simple adapted local descriptions of
geometrical objects with physical meaning
studied in Riemann-Lagrange geometry. At the same
time, these allow beautiful and natural adapted local descriptions of the
fundamental equations (generalized Maxwell or Einstein equations \cite{Ne8})
that dominate the abstract multi-time field theory created by
Riemann-Lagrange geometrical instruments.

Secondly, the use of adapted bases produced by {\it canonical nonlinear
connections} on 1-jet spaces
(i. e., which are built from the given multi-time dependent Lagrangian
function $L$ or the fundamental vertical metrical
d-tensor (\ref{fdt})) offers a {\it metrical character} to
Riemann-Lagrange theory of multi-time physical fields, according to
field theories classification from \cite{Got-Mars}.
In this direction, we point out that the construction of canonical
nonlinear connections on $ML^n_p$ spaces \cite{Ne5} is a {\it solved problem},
while the construction of canonical nonlinear connections on $GML^n_p$ spaces
\cite{Ne8} is still an {\it open problem}.
In what follows, we try to present a deep exposition of main geometrical results already
obtained, together with possible future directions of research, concerning
the construction of canonical nonlinear connections on $ML^n_p$ or
$GML^n_p$ spaces.

In this way, let us consider the first jet fibre bundle
$
J^1(T,M)\to T\times M
$, whose local coordinates
$(t^\alpha,x^i,x^i_\alpha)$, where $\alpha=\overline{1,p}$,
$i=\overline{1,n}$, transform by the rules

$$
\tilde t^\alpha=\tilde t^\alpha(t^\beta),
\quad
\tilde x^i=\tilde x^i(x^j),
\quad
\tilde x^i_\alpha={\partial\tilde x^i\over\partial x^j}{\partial t^\beta\over
\partial\tilde t^\alpha}x^j_\beta\;.
$$
\addtocounter{defin}{1}$\!\!$
{\bf Definition \thedefin} A pair $\Gamma=(M^{(i)}_{(\alpha)\beta}, N^{(i)}_
{(\alpha)j})$ consisting of local functions on 1-jet bundle $E=J^1(T,M)$,
whose transformation rules are given by
$$
\tilde M^{(j)}_{(\beta)\mu}{\partial\tilde t^\mu\over\partial
t^\alpha}=M^{(k)}_{(\gamma)\alpha}{\partial\tilde x^j\over\partial x^k}
{\partial t^\gamma\over\partial\tilde t^\beta}-{\partial\tilde x^j_\beta\over
\partial t^\alpha}\;,\qquad
\tilde N^{(j)}_{(\beta)k}{\partial\tilde x^k\over\partial
x^i}=N^{(k)}_{(\gamma)i}{\partial\tilde x^j\over\partial x^k}
{\partial t^\gamma\over\partial\tilde t^\beta}-{\partial\tilde x^j_\beta\over
\partial x^i}\;,
$$
is called a {\it nonlinear connection} on $E$. The components $M^{(i)}_{(\alpha)
\beta}$ (resp. $N^{(i)}_{(\alpha)j}$) are called
the {\it temporal} (resp.
{\it spatial}) {\it components} of the nonlinear connection $\Gamma$.
\medskip\\
\addtocounter{ex}{1}
{\bf Example \theex} If $h_{\alpha\beta}(t^\mu)$ (resp. $\varphi_{ij}(x^m)$)
is a semi-Riemannian metric on the temporal (resp. spatial) manifold $T$
(resp. $M$), and $H^\gamma_{\alpha\beta}(t^\mu)$ (resp. $\gamma^k_{ij}(x^m)$)
are its Christoffel symbols, then the pair of local functions
$
\Gamma_0=(\stackrel{0}{M}\;^{\!\!\!\!(j)}_{\!\!\!\!(\beta)\alpha},
\stackrel{0}{N}\;^{\!\!\!\!(j)}_{\!\!\!\!(\beta)i}),
$
where

\begin{equation}\label{can}
\stackrel{0}{M}\;^{\!\!\!\!(j)}_{\!\!\!\!(\beta)\alpha}=-H^\gamma_{\alpha
\beta}x^j_\gamma\;,
\qquad
\stackrel{0}{N}\;^{\!\!\!\!(j)}_{\!\!\!\!(\beta)i}=
\gamma^j_{ik}x^k_\beta\;,
\end{equation}
represents a nonlinear connection on $E$, which is called the
{\it canonical nonlinear connection of the semi-Riemannian metrics
$h_{\alpha\beta}$ and $\varphi_{ij}$}.\medskip

In this context, an important geometrical concept used in our studies, whose
physical meaning is intimately connected by the concept of energy, is
introduced by
\medskip\\
\addtocounter{defin}{1}
{\bf Definition \thedefin} A smooth map $f\in C^\infty(T,M)$, whose local
components verify the PDEs system of order two

$$
h^{\alpha\beta}\{x^i_{\alpha\beta}+M^{(i)}_{(\alpha)\beta}+N^{(i)}_
{(\alpha)m}x^m_\beta\}=0,
$$
where $h_{\alpha\beta}$ is a semi-Riemannian metric on $T$,
is called a {\it $h$-generalized harmonic map of the  nonlinear connection
$\Gamma=(M^{(i)}_{(\alpha)\beta}, N^{(i)}_{(\alpha)j})$}.
\medskip\\
\addtocounter{ex}{1}
{\bf Example \theex} The generalized harmonic maps of the nonlinear
connection (\ref{can}) are equivalent with classical harmonic maps between
the semi-Riemannian manifolds $(T,h)$ and $(M,\varphi)$. This statement
emphasizes the naturalness of previous definition.

\section{Canonical nonlinear connections on $ML^n_p$ spaces}

\setcounter{equation}{0}
\hspace{5mm} Concerning the construction of a canonical nonlinear connection
$\Gamma_L$ from a given multi-time dependent Lagrangian ${\cal L}=
L\sqrt{\vert h\vert}$,  let us consider the $ML^n_p$ space (i. e., the
metrical multi-time Lagrange space)
$ML^n_p=(J^1(T,M),L)$, where $\dim T=p$, $\dim M=n$, whose multi-time
dependent Lagrangian function $L$ is of the form \cite{Ne4}:

$$
L(t^\gamma,x^k,x^k_\gamma)=
\left\{\begin{array}{ll}\medskip
L(t,x^i,y^i),&p=1
\\
h^{\alpha\beta}(t^\gamma)g_{ij}(t^\gamma,x^k)x^i_\alpha x^j_\beta+U^{(\alpha)}
_{(i)}(t^\gamma,x^k)x^i_\alpha+F(t^\gamma,x^k),&p\geq 2,
\end{array}\right.
$$
where $U^{(\alpha)}_{(i)}(t^\gamma,x^k)$ is a d-tensor on $J^1(T,M)$ and $F$
is a smooth function on $T\times M$.

From our point of view, it is very interesting that the construction of
$\Gamma_L$ is strongly connected of the critical points of the energy
action functional of $L$, via its attached generalized harmonic maps.
To be more clearly, assume that the  semi-Riemannian temporal
manifold $(T,h)$ is compact and orientable. In this context, we define the
{\it multi-time
relativistic $h$-energy functional} of the Lagrangian function $L$:

$$
\E_L:C^\infty(T,M)\to\R,\quad
\E_L(f)=\int_TL(t^\alpha,x^i,x^i_\alpha)\sqrt{\vert h\vert}\;dt^1\wedge dt^2
\wedge\ldots\wedge dt^p,
$$
where the smooth map $f$ is locally expressed by $(t^\alpha)\to(x^i(t^\alpha))$.
Now, using an important result proved in \cite{Ne4}, we immediately find

\begin{th}\label{ml}
The extremals of the multi-time relativistic $h$-energy functional $\E_L$
of the metrical multi-time Lagrange space $ML^n_p$ are equivalent with the
$h$-generalized harmonic maps of the nonlinear connection $\Gamma_L$ defined by
the components:
\medskip\\
$$
\begin{array}{l}\bigskip
M^{(i)}_{(\alpha)\beta}=-H^\gamma_{\alpha\beta}x^i_\gamma,
\hspace*{85mm}
p\geq 1,
\\
N^{(i)}_{(\alpha)j}=
\left\{\begin{array}{ll}\medskip
\displaystyle{{g^{ik}\over 4}\left[{\partial^2L\over\partial x^j\partial
y^k}y^j-{\partial L\over\partial x^k}+{\partial^2L\over\partial t\partial y^k}+
{\partial L\over\partial x^k}H^1_{11}+2h^{11}H^1_{11}g_{kl}y^l\right]},&p=1
\\\medskip
\displaystyle{\Gamma^i_{jk}x^k_\alpha+{g^{ik}\over 2}{\partial g_{jk}\over
\partial t^\alpha}+{g^{ik}\over 4}h_{\alpha\gamma}U^{(\gamma)}_{(k)j}},&
p\geq 2,
\end{array}\right.
\end{array}
$$
where
$
\Gamma^l_{jk}=(g^{li}/2)(\partial g_{ij}/\partial x^k+\partial
g_{ik}/\partial x^j-\partial g_{jk}/\partial x^i)
$,
$
U^{(\alpha)}_{(i)j}=\partial U^{(\alpha)}_{(i)}/\partial x^j-
\partial U^{(\alpha)}_{(j)}/\partial x^i.
$
\end{th}
\addtocounter{defin}{1}
{\bf Definition \thedefin}
The nonlinear connection $\Gamma_L$ from Theorem \ref{ml}
is called the {\it canonical nonlinear connection of the metrical multi-time Lagrange
space $ML^n_p$.}
\medskip\\
\addtocounter{rem}{1}
{\bf Remarks \therem} i) In the particular case of usual time axis
$(T,h)=(\R,\delta)$, the canonical nonlinear connection
$\Gamma_L=(0,N^{(i)}_{(1)j})$, produced by the time dependent
Lagrangian ${\cal L}=L\sqrt{\vert h\vert}$ of a {\it relativistic
rheonomic Lagrange space} \cite{Ne6}

$$
RL^n=(J^1(\R,M)\equiv\R\times\T M,L),
$$
generalizes, in relativistic dynamical terms, canonical nonlinear
connections used in theory of classical rheonomic Lagrange spaces
\cite{Mir-An}, Ch. XIII.

ii) Note that, in the case $p\geq 2$, the construction of $\Gamma_L$ on a
$ML^n_p$ space essentially relies on Kronecker
$h$-regularity of the fundamental vertical metrical d-tensor produced by
$L$ (i. e.,
$
G^{(\alpha)(\beta)}_{(i)(j)}=(1/2)(\partial^2L/\partial x^i_\alpha\partial
x^j_\beta)=h^{\alpha\beta}(t^\gamma)g_{ij}(t^\gamma,x^k))
$.

\section{Canonical nonlinear connections on $GML^n_p$ spaces}

\setcounter{equation}{0}
\hspace{5mm} More general, in order to construct a canonical nonlinear
connection $\Gamma_G$ on a generalized metrical multi-time Lagrange space
$GML^n_p=(J^1(T,M),G^{(\alpha)(\beta)}_{(i)(j)})$, whose Kronecker
$h$-regular vertical fundamental metrical d-tensor (\ref{fdt})
is not necessarily provided by a multi-time dependent Lagrangian function $L$,
we point out that the temporal semi-Riemannian metric $h_{\alpha\beta}$
naturally produces the {\it temporal components} of the nonlinear connection
$\Gamma_G$, taking

\begin{equation}\label{tc}
M^{(i)}_{(\alpha)\beta}=\stackrel{0}{M}\;^{\!\!\!\!(i)}_{\!\!\!\!(\alpha)
\beta}=-H^\gamma_{\alpha\beta}x^i_\gamma.
\end{equation}

Concerning the construction of the {\it spatial components}
of the canonical nonlinear connection $\Gamma_G$,
we emphasize that the Riemann-Lagrange geometrical background
for the generalized multi-time field theory from \cite{Ne8} relies on the
use of an {\it "a priori"} given {\it without torsion} nonlinear connection
$\Gamma$, that is the spatial components of $\Gamma$ verify the equalities

$$
{\partial N^{(i)}_{(\alpha)j}\over\partial x^k_\gamma}=
{\partial N^{(i)}_{(\alpha)k}\over\partial x^j_\gamma}.
$$
Consequently, in order to construct the canonical nonlinear connection
$\Gamma_G$ of a $GML^n_p$ space, we need to produce some without torsion
spatial components from the fundamental vertical metrical d-tensor
(\ref{fdt}). In this direction, following geometrical ideas from previous
section, we introduce the concept of {\it multi-time relativistic
energy  Lagrangian function} of a $GML^n_p$ space, setting

$$
E_G=G^{(\mu)(\nu)}_{(m)(r)}x^m_\mu x^r_\nu=h^{\mu\nu}(t^\gamma)g_{mr}(t^\gamma,x^k,x^k_
\gamma)x^m_\mu x^r_\nu\;,
$$
together with its {\it multi-time relativistic
$\psi$-energy action functional}

$$
\E_G:C^\infty(T,M)\to\R,\quad
\E_G(f)=\int_T E_G(t^\gamma,x^k(t^\gamma),x^k_\gamma(t^\mu))
\sqrt{\vert\psi\vert}dt,
$$
where $\psi$ is a semi-Riemannian metric on the temporal manifold $T$.

Firstly, suppose that $E_G$ is a Kronecker $\psi$-regular
Lagrangian function of the form
$$
E_G=h^{\alpha\beta}(t^\gamma)g_{ij}(t^\gamma,x^k,x^k_
\gamma)x^i_\alpha x^j_\beta=\psi^{\alpha\beta}(t^\gamma)\varepsilon_{ij}
(t^\gamma,x^k)x^i_\alpha x^j_\beta\;,
$$
where $\varepsilon_{ij}(t^\gamma,x^k)$ is a symmetric d-tensor of rank $n$,
having a constant signature on $T\times M$.
In this context, applying Theorem \ref{ml} to the particular $ML^n_p$ space
$MLGML^n_p=(J^1(T,M),E_G)$ produced by the generalized metrical multi-time
Lagrange space $GML^n_p$, we can construct certain without torsion spatial
components
for $\Gamma_G$, via the $\psi$-energy functional $\E_G$. Consequently, we may
consider that the spatial components of the canonical nonlinear connection
$\Gamma_G$ are given by the formulas:

\begin{equation}\label{sc}
N^{(i)}_{(\alpha)j}=\Psi^i_{jm}x^m_\alpha+{\varepsilon^{im}\over 2}{\partial
\varepsilon_{jm}\over\partial t^\alpha},
\end{equation}
where
$
\Psi^i_{jk}(t^\mu,x^m)
$
are {\it generalized Christoffel symbols} for $\varepsilon_{ij}(t^\gamma,x^k)$.
\medskip\\
\addtocounter{rem}{1}
{\bf Remark \therem} A canonical nonlinear connection $\Gamma_L$,
derived from a multi-time dependent quadratic Lagrangian function of general
form

$$
L=G^{(\alpha)(\beta)}_{(i)(j)}(t^\gamma,x^k)x^i_\alpha x^j_\beta+U^{(\alpha)}
_{(i)}(t^\gamma,x^k)x^i_\alpha+F(t^\gamma,x^k)
$$
and a given temporal semi-Riemannian metric $h_{\alpha\beta}(t^\gamma)$, may
be naturally produced, using preceding ideas for the metrical d-tensor
$
{\cal G}^{(\alpha)(\beta)}_{(i)(j)}=h^{\alpha\beta}h_{\mu\nu}G^{(\mu)(\nu)}_
{(i)(j)}
$.\medskip

Secondly, suppose that $E_G$ is not a Kronecker $\psi$-regular
Lagrangian function for any semi-Riemannian metric $\psi$. Under these
assumptions, we are forced to give {\it "ab initio"} the without torsion
spatial components of the nonlinear connection $\Gamma_G$.\medskip\\
\addtocounter{ex}{1}
{\bf Example \theex} For the study of particular $GML^n_p$ spaces
$GRGML^n_p$ (i. e.,
that represents a geometrical model for multi-time Relativity and
Electromagnetism \cite{Ne9}) and $RGOGML^n_p$ (i. e., that represents a
geometrical background for dynamical relativistic multi-time optics
\cite{Ne10}), it is convenient to use the {\it "a priori"} given
without torsion spatial components

$$
N^{(i)}_{(\alpha)j}=\stackrel{0}{N}\;^{\!\!\!\!(i)}_{\!\!\!\!(\alpha)j}=
\gamma^i_{jm}x^m_\alpha
$$
which are directly provided by the Kronecker $h$-regular fundamental
vertical metrical d-tensors of
these spaces. The main physical motives, together with their intrinsec
difficulties, that determined us to use these spatial components in the
geometrical study of the spaces $GRGML^n_p$ and $RGOGML^n_p$,
are deeply exposed in \cite{Ne9}, \cite{Ne10}.\medskip

Taking  into account the preceding discussions related to the construction of
canonical nonlinear connections on $GML^n_p$ spaces, we consider important
to study the conditions that must be imposed to the multi-time relativistic
energy Lagrangian function $E_G$ of a generalized metrical multi-time Lagrange
space, in order to obtain its Kronecker $\psi$-regularity, where $\psi$ is an
arbitrary semi-Riemannian metric on $T$. We recall that, under the
Kronecker $\psi$-regularity assumption, $E_G$ produces without torsion
spatial components for the canonical nonlinear connection $\Gamma_G$ of the space
$GML^n_p$, via the formulas (\ref{sc}).

In order to do this study, let us suppose $p=\dim T\geq  2$. In this case,
the characterization theorem of $ML^n_p$ spaces from \cite{Ne4} imply
that the Kronecker $\psi$-regularity of $E_G$ reduces to the
existence of certain d-tensors $\varepsilon_{ij}(t^\gamma,x^k)$, $U^{(\alpha)}_
{(i)}(t^\gamma,x^k)$ on $J^1(T,M)$, together with a smooth function
$F:T\times M\to\R$, such that

$$
E_G=\psi^{\alpha\beta}(t^\gamma)\varepsilon_{ij}(t^\gamma,x^k)x^i_\alpha
x^j_\beta+U^{(\alpha)}_{(i)}(t^\gamma,x^k)x^i_\alpha+F(t^\gamma,x^k).
$$
In other words, $E_G$ is a Kronecker $\psi$-regular Lagrangian function
if the following equality is true for some $\varepsilon_{ij}$, $U^{(\alpha)}_
{(i)}$ and $F$:

$$
h^{\alpha\beta}(t^\gamma)g_{ij}(t^\gamma,x^k,x^k_\gamma)x^i_\alpha
x^j_\beta=\psi^{\alpha\beta}(t^\gamma)\varepsilon_{ij}(t^\gamma,x^k)x^i_\alpha
x^j_\beta+U^{(\alpha)}_{(i)}(t^\gamma,x^k)x^i_\alpha+F(t^\gamma,x^k).
$$

Suppose that the spatial metrical d-tensor $g_{ij}$ of the
space $GML^n_p$ does not depend on partial directions  $x^i_\alpha$.
It is obvious that,  taking
$\psi^{\alpha\beta}=h^{\alpha\beta}$, $\varepsilon_{ij}=g_{ij}$,
$U^{(\alpha)}_{(i)}=0$ and $F=0$, we can conclude that $E_G$ is a Kronecker
$h$-regular Lagrangian function. Therefore, in this situation, using the
formulas (\ref{sc}), we can build the spatial components of the canonical
nonlinear connection $\Gamma_G$.
\medskip\\
{\bf Open problems.} i) In the case $p=\dim T\geq 2$, are there
spatial metrical d-tensors $g_{ij}$ on $GML^n_p$ spaces, depending effectively
on partial directions $x^i_\alpha$, such that $E_G$ to be a Kronecker
$\psi$-regular Lagrangian function?

ii) Is it possible  to construct, in a natural way, without torsion spatial
components of a nonlinear connection $\Gamma_G$, which to be canonically
produced by a Kronecker $h$-regular vertical fundamental metrical d-tensor
$G$ on $J^1(T,M)$?

{\it Author's address:}
M. Neagu,
University Politehnica of Bucharest,
Department of Mathematics I,
Splaiul Independentei 313,
77206 Bucharest, Romania.\\
{\it E-mail:} mircea@mathem.pub.ro

\end{document}